\documentclass[12pt]{amsart}
\usepackage{amssymb,amsmath,amsthm,fullpage}
\usepackage{pdfsync,enumerate,color}

\newcommand{\abs}[1]{\left\vert#1\right\vert}

\renewcommand{\L}{\mathcal{L}}

\newcommand{\R}{\mathbb R}

\newcommand{\C}{\mathbb C}

\DeclareMathOperator{\Imm}{Im}

\DeclareMathOperator{\Rre}{Re}
\DeclareMathOperator{\Dom}{Dom}

\newcommand{\bd}{\textrm{b}}

\newcommand{\p}{\partial}

\newcommand{\z}{\bar z}
\newcommand{\w}{\bar w}
\newcommand{\dbar}{\bar\partial}
\newcommand{\dbars}{\bar\partial^*}
\newcommand{\dbarst}{\bar\partial^*_t}

\DeclareMathOperator{\Tr}{Tr}

\newcommand{\la}{\langle}
\newcommand{\ra}{\rangle}

\newcommand{\LL}{\bar L}
\DeclareMathOperator{\dist}{dist}

\newtheorem{thm}{Theorem}[section]

\newtheorem{lem}[thm]{Lemma}
\newtheorem{cor}[thm]{Corollary}

\newtheorem*{theorem*}{Theorem}

\theoremstyle{definition}
\newtheorem{defn}[thm]{Definition}

\theoremstyle{remark}
\newtheorem{rem}[thm]{Remark}

\newcommand{\Om}{\Omega}

\begin{document}

\title[Boundary estimates on $Z(q)$ domains]{A remark on boundary estimates on unbounded $Z(q)$ domains in $\C^n$}%

\author{Phillip S. Harrington and Andrew Raich}%
\address{SCEN 309, 1 University of Arkansas, Fayetteville, AR 72701}%
\email{psharrin@uark.edu, araich@uark.edu}%

\thanks{The second author is partially supported by NSF grant DMS-1405100.}%

\begin{abstract}The goal of this note is to explore the relationship between the Folland-Kohn basic estimate and the $Z(q)$-condition. In particular, on unbounded pseudoconvex (resp., pseudoconcave) domains,
we prove that the Folland-Kohn basic estimate is equivalent to uniform strict pseudoconvexity (resp., pseudoconcavity).
As a corollary, we observe that despite the Siegel upper half space being strictly pseudoconvex and biholomorphic to a the unit ball,
it fails to satisfy uniform strict pseudoconvexity and hence the Folland-Kohn basic estimate fails.

On unbounded non-pseudoconvex domains, we show that the Folland-Kohn basic estimate on $(0,q)$-forms
implies a uniform $Z(q)$ condition, and conversely, that a uniform $Z(q)$ condition with some additional hypotheses implies the Folland-Kohn basic estimate for $(0,q)$-forms.
\end{abstract}

\maketitle

\section*{Acknowledgement}
We would like to dedicate this note to our colleague John Ryan on the occasion of his sixtieth birthday.

%%%%%%%%%%%%%%%%%%%
%
%				SECTION: INTRODUCTION
%
%%%%%%%%%%%%%%%%%%%%
\section{Introduction} The basic estimate in the sense of Folland and Kohn is a bound of the weighted boundary $L^2$ norm of a $(0,q)$-form $f$ by the
weighted $L^2$ norms of $f$, $\dbar f$, and $\dbars_t f$.
In this paper, we prove that on an unbounded pseudoconvex (resp., pseudoconcave) domain $\Omega\subset\C^n$, the basic estimate holds in the sense of Folland and Kohn \cite{FoKo72} if and only if
$\Omega$ is uniformly pseudoconvex (resp., pseudoconcave). In the non-pseudoconvex case, we show that if $\Omega$ satisfies the Folland-Kohn basic estimate
on $(0,q)$-forms, then $\Omega$ satisfies a
uniform $Z(q)$ condition. We also establish a partial converse. If $\Omega$ satisfies a uniform $Z(q)$ condition, together with the hypotheses that guarantee the closed range estimate from
\cite{HaRa15u}, then $\Om$ satisfies a strong form of the Folland-Kohn basic estimate for $(0,q)$-forms.

In Theorems 3.2.1 and 3.2.4 in \cite{Hor65}, H\"ormander showed that condition $Z(q)$ is locally equivalent to the basic estimate on $(0,q)$-forms in the sense of Folland and Kohn \cite{FoKo72}, with a global equivalence on bounded domains (see \cite[Corollary 3.2.3]{Hor65}).  The motivation for this paper is to show that this equivalence fails on unbounded domains but H\"ormander's Theorem 3.2.1 actually implies a stronger condition than $Z(q)$ on unbounded domains.  Specifically, the
$Z(q)$ condition must be replaced with a uniform $Z(q)$ condition. These two conditions are equivalent on bounded domains.

When $q=1$, the uniform $Z(q)$ condition is a uniform strict pseudoconvexity condition which we will show that the Siegel upper half space fails to satisfy. As a corollary, we
obtain a result that is strikingly different than the bounded case, namely, strict pseudoconvexity is no longer sufficient for the Folland-Kohn basic estimate to hold.

This note is continuation of our exploration of the $L^2$ theory of $\dbar$ on unbounded domains. The $L^2$-theory of the $\dbar$ operator in $L^2_{0,q}(\Omega)$ is defective for most unbounded domains -- see \cite{HaRa15u,HeMcN15}. This fact motivated
the development of a function theory that could overcome the problems with the unweighted theory.
In  \cite{HaRa13,HaRa14}, we built a theory for Sobolev spaces and elliptic operators that
is appropriate for an examination of the $\dbar$-Neumann problem on unbounded domains. In \cite{HaRa15u}, we formulated a weak $Z(q)$ condition that sufficed to prove closed range for
$\dbar$ and (consequently) the continuity of a weighted $\dbar$-Neumann operator $N_{q,t}$ on $(0,q)$-forms in an appropriate weighted Sobolev space $H^s_t(\Omega)$.

Gansberger was the first to investigate closed range of $\dbar$ in (weighted) $L^2$ on unbounded domains \cite{Gan10}, but his focus was on compactness of the $\dbar$-Neumann operator
on pseudoconvex domains. In \cite{HaRa11}, we began our investigation of sufficient conditions for closed range of $\dbar$ on $(0,q)$-forms for a \emph{fixed} $q$, $0<q<n$. This led to a generalization
in \cite{HaRa15}, suitable for $C^3$ domains in a Stein manifold. It was the definition of weak $Z(q)$ in \cite{HaRa15} that we modified for unbounded domains in \cite{HaRa14}.
Herbig and McNeal \cite{HeMcN15} also explore closed range of $\dbar$ on pseudoconvex domains that are not necessarily bounded.

The outline of this note is as follows: we state definitions and formulate the main results in Section \ref{sec:main results} and prove the main theorems in Section \ref{sec:applications}.

%%%%%%%%%%%%%%%%%%%%%%
%
%				SECTION: MAIN RESULTS
%
%%%%%%%%%%%%%%%%%%%%%%
\section{Definitions and Results}\label{sec:main results}

% subsection: Definitions
\subsection{Definitions}\label{subsec:definitions}
To continue the discussion, we need to introduce some terminology. Our definitions follow the setup in \cite{HaRa13,HaRa14,HaRa15u}.

A function $\rho:\C^n\to\R$ is called a \emph{defining function} for $\Omega$ if $\Omega = \{z : \rho(z)<0\}$ and
$d\rho\neq 0$ on $\bd\Omega$. The \emph{Levi form} of $\Omega$
is the restriction of the complex Hessian of a defining function $\rho$ to the maximal complex tangent space (locally, an $(n-1)\times(n-1)$ matrix).  The induced CR-structure on $\bd\Omega$  at $z\in\bd\Omega$ is
\[
T^{1,0}_z(\bd\Omega)  = \{ L\in T^{1,0}(\C) : \p\rho(L)=0 \}.
\]
Let $T^{1,0}(\bd\Omega)$ be the space of $C^{m-1}$ sections of $T^{1,0}_z(\bd\Omega)$ and $T^{0,1}(\bd\Omega) = \overline{T^{1,0}(\bd\Omega)}$.
We denote the exterior algebra generated by the dual forms by $\Lambda^{p,q}(\bd\Omega)$.
If we normalize $\rho$ so that $|d\rho|=1$ on $\bd\Omega$, then the \emph{normalized Levi form} $\L$ is the real element of $\Lambda^{1,1}(\bd\Omega)$ defined by
\[
\L(-i L\wedge \LL) = i\p\dbar\rho(-iL\wedge\LL)
\]
for any $L\in T^{1,0}(\bd\Omega)$.

% definition: uniform Z(q)
\begin{defn}
Let $\Omega\subset\C^n$ be a domain with a $C^m$ boundary, $m\geq 2$. We say that $\Omega$ satisfies
\emph{condition $Z(q)$} if the Levi form has at least $(n-q)$-positive or at least $(q+1)$-negative eigenvalues. We say that $\Omega$ satisfies the
\emph{uniform $Z(q)$ condition} if for some $\lambda>0$ the normalized Levi form has at least $(n-q)$ eigenvalues greater than $\lambda$ or at least $(q+1)$ eigenvalues smaller than $-\lambda$.
\end{defn}

Let $\Omega\subset \C^n$ be a domain with $C^m$ boundary $\bd\Omega$, $m\geq 3$.
\begin{defn}\label{defn:uniform_defining_function}
  Let $\Omega\subset\mathbb{R}^n$, and let $\rho$ be a $C^m$ defining function for $\Omega$ defined on a neighborhood $U$ of $b\Omega$ such that
  \begin{enumerate}
    \item $\dist(\partial\Omega,\partial U)>0$,

    \item $\|\rho\|_{C^m(U)}<\infty$,

    \item $\inf_U|\nabla\rho|>0$.
  \end{enumerate}
  We say that such a defining function is \emph{uniformly $C^m$}.  If $\rho$ on $U$ is uniformly $C^m$ for all $m\in\mathbb{N}$, we say $\rho$ is \emph{uniformly $C^\infty$}.
\end{defn}
In \cite{HaRa13}, we show that we may assume $\abs{\nabla\rho}=1$ on $U$ without loss of generality.  In fact, the existence of any uniformly $C^m$ defining function implies that the signed distance function is uniformly $C^m$.

We identify real $(1,1)$-forms with a hermitian matrix as follows:
\[
c=\sum_{j,k=1}^{n} i c_{j\bar k}\, dz_j\wedge d\z_k
\]

%For a function $\alpha$, we denote $\alpha_k = \frac{\p\alpha}{\p z_k}$ and $\alpha_{\bar j} = \frac{\p\alpha}{\p\z_j}$.
%Let
%\[
%\Omega_{\ep,\rho} = \{z\in \Omega : -\ep <\rho(z)< 0\}
%\qquad\text{and}\qquad
%\Omega_{\ep,\rho}' = \{z\in \C^n : -\ep < \rho(z) < \ep\}.
%\]
We denote the $L^2$-inner product on $L^2(\Omega,e^{-t|z|^2})$ by
\[
(f,g)_t = \int_\Omega \la f, g\ra\, e^{-t|z|^2} dV = \int_\Omega f \bar g\, e^{-t|z|^2} dV
\]
where $\la \cdot,\cdot \ra$ is the standard pointwise inner product on $\C^n$ and $dV$ is Lebesgue measure on $\C^n$. We denote
$d\sigma$ as the induced surface area measure on $\bd\Omega$. Also $\| f \|_t^2 = \int_\Omega |f|^2 e^{-t|z|^2}\, dV$.

Let $\dbar:L^2_{0,q}(\Omega, e^{-t|z|^2}) \to L^2_{0,q+1}(\Omega, e^{-t|z|^2})$ denote the maximal closure of the Cauchy-Riemann operator, and let $\dbars_t:L^2_{0,q+1}(\Omega, e^{-t|z|^2}) \to L^2_{0,q}(\Omega, e^{-t|z|^2})$ denote the adjoint with respect to $(\cdot,\cdot)_t$.

%Let $\I_q = \{ (i_1,\dots,i_q)\in \N^n : 1 \leq i_1 < \cdots < i_q\leq n \}$. For $I\in\I_{q-1}$, $J\in\I_q$, and $1\leq j \leq n$, let
%$\ep^{jI}_{J} = (-1)^{|\sigma|}$ if $\{j\} \cup I = J$ as sets and $\sigma$ is the length of the permutation that takes $\{j\}\cup I$ to $J$. Set
%$\ep^{jI}_J=0$ otherwise. We use the standard notation that if $u = \sum_{J\in\I_q} u_J\, d\z_J$, then
%\[
%u_{jI} = \sum_{J\in\I_q} \ep^{jI}_J u_J.
%\]

%Let $L^{t}_j = \frac{\p}{\p z_j} - t\z_j = e^{t|z|^2}\frac{\p}{\p z_j} e^{-t|z|^2}$ and let $\dbars_t:L^2_{0,q+1}(\Omega, e^{-t|z|^2}) \to L^2_{0,q}(\Omega, e^{-t|z|^2})$ be the $L^2$-adjoint of
%$\dbar:L^2_{0,q}(\Omega, e^{-t|z|^2}) \to L^2_{0,q+1}(\Omega, e^{-t|z|^2})$. This means that if
%$f = \sum_{J\in\I_q} f_J\, d\z_J$ and $g = \sum_{L\in\I_{q+1}}g_L\, d\z_L\in\Dom(\dbars_t)$, then
%\[
%\dbar f = \sum_{\atopp{J\in\I_q}{L\in\I_{q+1}}}\sum_{k=1}^n \ep^{kJ}_L \frac{\p f_J}{\p\z_k}\, d\z_L
%\qquad\text{and}\qquad
%\dbars_t g = -\sum_{J\in\I_{q}} \sum_{j=1}^n L^t_j f_{jJ}\, d\z_J.
%\]

% Definition: the Folland-Kohn basic estimate
\begin{defn}\label{defn:the FK basic est}
Let $\Omega\subset\C^n$ be a domain of class $C^1$ and $1 \leq q \leq n$. We say that $\Omega$ satisfies the \emph{Folland-Kohn basic estimate on $(0,q)$-forms} if
there exists $t>0$ and a constant $C_t>0$ so that for all $(0,q)$-forms $f\in \Dom(\dbar)\cap\Dom(\dbarst)$, $1 \leq q \leq n-1$,
\begin{equation}\label{eqn:bdy dominated by dbar,dbars}
\int_{\bd\Omega} |f|^2 e^{-t|z|^2}\, d\sigma \leq C_t \Big(\| \dbar f\|_t^2 + \| \dbars_t f \|_t^2\Big).
\end{equation}
\end{defn}

%If $U$ is a suitably small neighborhood of $\bd\Omega$, we
%use $\tau$ to denote the orthogonal projection and
%restriction
%\[
%\tau : \Lambda^{p,q}(U)\to \Lambda^{p,q}(\bd\Omega).
%\]

% Definition: tubular neighborhood
%\begin{defn}\label{defn:tubular nbhd}Given a set $M\subset\C^n$, a \emph{tubular neighborhood} of $M$ is an open set $U_r$ of the form
%$U_{r} = \{p\in\C^n : \dist(p,M)<r\}$ where $\dist(\cdot,\cdot)$ is the Euclidean distance function. We call $r$ the \emph{radius} of $U_r$.
%\end{defn}

Our definition for weak $Z(q)$ follows \cite{HaRa15}.
% Definition: weak Z(q)
\begin{defn}\label{defn:weak Z(q)}
Let $\Omega\subset \C^n$ be a domain with a uniformly $C^m$ defining function $\rho$ satisfying $|d\rho||_{b\Omega}=1$, $m\geq 3$. We say $\bd\Omega$ (or $\Omega$) satisfies
\emph{Z(q) weakly}
if there exists a hermitian matrix $\Upsilon=(\Upsilon^{\bar k j})$ of functions on $b\Omega$ that are uniformly bounded in $C^{m-1}$ such that %$\sum_{j=1}^{n}\Upsilon^{\bar k j}\rho_j=0$ on $b\Omega$ and:
\begin{enumerate}\renewcommand{\labelenumi}{(\roman{enumi})}
 \item All eigenvalues of $\Upsilon$ lie in the interval $[0,1]$.

 \item $\mu_1+\cdots+\mu_q-\sum_{j,k=1}^n\Upsilon^{\bar k j}\rho_{j\bar k}\geq 0$ where  $\mu_1,\ldots,\mu_{n-1}$ are the eigenvalues of the Levi form $\L$ in increasing order.

 \item $ \inf_{z\in\bd\Omega} \{ |q-\Tr(\Upsilon)|\} >0$.

\end{enumerate}
%We say that $\Omega$ satisfies \emph{Z(q)} if the inequality in (ii) is strictly positive and \emph{Z(q) uniformly} if the inequality in (ii) is bounded from below a constant $\mu_{\min}>0$.
\end{defn}
%\begin{remark}
%The hypothesis on $\Omega$ ensures that there exists an $r>0$ so that
%if $p\in\bd\Omega$, then $B(p,r)\cap \bd\Omega$ is connected. This guarantees than we can move away from the boundary a uniform distance
%without intersecting another piece of $\bd\Omega$.
%\end{remark}

We showed in \cite{HaRa15u} that weak $Z(q)$ suffices for the closed range of $\dbar$ in $L^2_{0,q}(\Omega,e^{-t|z|^2})$ and $L^2_{0,q+1}(\Omega,e^{-t|z|^2})$ for
$t$ (when the Levi-form has $n-q$ nonnegative eigenvalues) or $-t$ (when the Levi-form has $q+1$ nonpositive eigenvalues) sufficiently large.
Additionally, we proved the following lemma in \cite{HaRa15} outlining the relationship between weak $Z(q)$ and $Z(q)$.
\begin{lem}[Lemma 2.8, \cite{HaRa15}]
  \label{lem:count_eigenvalues}
  For $1\leq q\leq n-1$ let $\Omega\subset M$ be a domain for which $b\Omega$ satisfies $Z(q)$ weakly.  Let $\Upsilon$ be as in Definition \ref{defn:weak Z(q)}.
   For any fixed boundary point, if $\Tr(\Upsilon)<q$ then the Levi form has at least $n-q$ nonnegative eigenvalues, and if $\Tr(\Upsilon)>q$, then the Levi form has at least $q+1$ nonpositive eigenvalues.

   Additionally, if the inequality in Definition \ref{defn:weak Z(q)}.ii  is strictly positive, then $\Om$ satisfies $Z(q)$, and if this inequality is uniformly bounded away from $0$, then $\Om$ satisfies  $Z(q)$ uniformly.
\end{lem}
In \cite{HaRa15}, Lemma \ref{lem:count_eigenvalues} is stated for boundary terms but the result is local and boundedness appears nowhere in the proof. Also, the statements in the lemma about
$Z(q)$ and uniform $Z(q)$ are implicit in the argument leading to \cite[Lemma 2.8]{HaRa15}.
%In \cite{HaRa14}, we introduced six hypotheses $(HI)-(HVI)$ that were important for developing elliptic theory with weighted Sobolev spaces on unbounded domains.  The first hypothesis was equivalent to Definition \ref{defn:uniform_defining_function}, so $(HI)$ will be satisfied whenever we have a uniformly $C^m$ defining function with $m\geq 3$.  Hypotheses $(HII)-(HV)$ are trivial for the weight function $\varphi=|z|^2$, so we will not need to address them directly in this paper.  Thus, we need only concern ourselves with $(HVI)$.  In the notation of the present paper, we have:

%To use the Sobolev space machinery of \cite{HaRa14} with the weight $e^{-t|z|^2}$, we need only the following condition at $\infty$.
%\begin{defn}
%  Let $\Omega\subset\mathbb{R}^{2n}$ be an unbounded domain.  We say $\Omega$ is \emph{asymptotically non-radial} if
%  \[
%    \inf_{r>0}\sup_{|x|>r,x\in\bd\Omega}\frac{x\cdot\nabla\rho}{|x||\nabla\rho|}<1
%  \]
%  for any $C^1$ defining function $\rho$ for $\Omega$.
%\end{defn}

% subsection: Results
\subsection{Results}\label{subsec:results}
% Theorem: uniform Z(q) equivalent to controlling boundary norm
\begin{thm}\label{thm:surface int, q-convex}
Let $\Omega\subset\C^n$ be a domain of class $C^3$ such that the sum of any $q$ eigenvalues of the normalized Levi-form is nonnegative or the sum of any $n-1-q$ eigenvalues of the normalized Levi-form is nonpositive for some $1\leq q\leq n-1$.  The domain $\Omega$
satisfies the uniform $Z(q)$ condition if and only if there exists $t_0\geq 0$ such that $\Omega$ satisfies the Folland-Kohn basic estimate on $(0,q)$-forms for all $t\geq t_0$.  When the sum of any $q$ eigenvalues of the normalized Levi-form is nonnegative, we only require $C^2$ boundary and we may take $t_0=0$.
\end{thm}

\begin{rem}
  The assumption that the sum of any $q$ eigenvalues of the normalized Levi-form is nonnegative is known as weak $q$-convexity in \cite{Ho91}.
\end{rem}

\begin{cor}\label{cor:surface int, pseudoconvex}
Let $\Omega\subset\C^n$ be a pseudoconvex (resp., pseudoconcave) domain of class $C^3$. The domain $\Omega$
is uniformly strictly pseudoconvex (resp., pseudoconcave) if and only if there exists $t_0\geq 0$ such that $\Omega$ satisfies the Folland-Kohn basic estimate on $(0,q)$-forms for all $1 \leq q \leq n-1$ and $t\geq t_0$.  When $\Omega$ is pseudoconvex, we only require $C^2$ boundary and we may take $t_0=0$.
\end{cor}

To specialize to $(0,q)$-forms, we need to appeal to the full basic estimate in \cite{HaRa15u}.
% Theorem: uniform Z(q) equivalent to controlling boundary norm
\begin{thm}\label{thm:surface int equiv to Z(q)}
Let $\Omega\subset\C^n$ be a $C^2$ domain and $1 \leq q \leq n-1$. If $\Omega$ satisfies the Folland-Kohn basic estimate for $(0,q)$-forms, then $\Omega$ satisfies
$Z(q)$ uniformly.

Conversely, if $\Omega$ is an
an unbounded domain admitting a uniformly $C^3$ defining function and satisfies  $Z(q)$ uniformly and $Z(q)$ weakly, then there exists $t_0\geq 0$ so that
whenever $t>t_0$, there exists $C_t$ so that $\Omega$ satisfies the Folland-Kohn basic estimate on $(0,q)$-forms for all $t\geq t_0$.
\end{thm}

Our primary (non)example is the Siegel upper half space
\[
\Sigma = \{(z,w)\in \C^{n+1} : \Imm w > |z|^2\}.
\]
The boundary of the Siegel upper half space is the Heisenberg group, $\mathbb H^n$.
Given the description of $\Sigma$ as a global graph, it is natural to use the defining function $r(z,w) = |z|^2 - \frac{w-\w}{2i}$, but $r$ is not a uniformly $C^3$ defining function since
$|dr| = 1 + 4|z|^2$ is not bounded.  However, $\rho(z,w)=\frac{r(z,w)}{|dr(z,w)|}$ is a uniformly $C^\infty$ defining function.
% Theorem: Siegel upper half space not uniform Z(q)
\begin{thm}\label{thm:Siegel upper half space}
The Siegel upper half space $\Sigma$ is a uniformly $C^\infty$ strictly pseudoconvex domain, but is not uniformly strictly pseudoconvex, so the Folland-Kohn basic
estimate (\ref{eqn:bdy dominated by dbar,dbars}) fails to hold on $(0,q)$-forms for any $1 \leq q \leq n-1$ and $t\geq 0$.
\end{thm}

\section{Proofs and Examples}\label{sec:applications}
In this section, we prove Theorem \ref{thm:Siegel upper half space} and provide an example of a uniformly $Z(q)$ domain. The proofs of Theorem \ref{thm:surface int, q-convex}
and Theorem \ref{thm:surface int equiv to Z(q)} are essentially the same -- the main difference is which basic estimate from \cite{HaRa15u} to apply.

\subsection{Proofs of Theorem \ref{thm:surface int, q-convex} and Theorem \ref{thm:surface int equiv to Z(q)}}
The proof that \eqref{eqn:bdy dominated by dbar,dbars} implies the uniform $Z(q)$ condition can be taken nearly verbatim from the argument that establishes \cite[Corollary 3.2.3]{Hor65}.
H\"ormander's arguments reduce to a local argument that renders irrelevant both the weight and unboundedness.

The proofs of the converse directions follow from the basic estimates in \cite{HaRa15u}. For Theorem \ref{thm:surface int, q-convex} in the case where the sum of any $q$ eigenvalues of the Levi-form are nonnegative, use \cite[Proposition 3.1]{HaRa15u} and \cite[Lemma 4.7]{Str10}.  When the sum of any $n-q-1$ eigenvalues is nonpositive, use \cite[Proposition 3.3]{HaRa15u} and \cite[Lemma 3.4]{HaRa15u} with $\Upsilon$ equal to the identity restricted to $T^{1,0}(b\Omega)$.  For Theorem \ref{thm:surface int equiv to Z(q)}, use \cite[Proposition 3.3]{HaRa15u} and \cite[Lemma 3.4]{HaRa15u}. The only difference for this proof is that our hypotheses are sufficient
to force the boundary integral term to dominate $\int_{\bd\Om} |f|^2\, e^{-t|z|^2}\, dz$ using \cite[Lemma 4.7]{Str10}. When $\Om$ is pseudoconcave, we have a basic identity analogous to \cite[Proposition 3.1]{HaRa15u} but given by
\cite[Equation (2)]{HaPeRa15}, with the more general basic identity given at the end of \cite[Section 1.2]{HaPeRa15}.
\qed

% subsection: Proof of Theorem \ref{thm:Siegel upper half space}
\subsection{Proof of Theorem \ref{thm:Siegel upper half space}\label{subsec:Heisenberg group proof}}
The definition of the Levi form required $|d\rho|=1$ on $\bd\Sigma$. The standard defining function for the Siegel upper half space $\Sigma = \{(z,w):\Imm w = |z|^2 \}$ is
\[
r(z,w) = \frac{w-\w}{2i}-|z|^2.
\]
This function is clearly unsuitable to use to compute the normalized Levi form. Instead, we use the signed distance function, $\tilde\delta(z)$,  defined by
\[
\tilde\delta(z) = \begin{cases} \dist(z,\bd\Sigma) &\text{if }z\not\in\bar\Sigma \\ -\dist(z,\bd\Sigma) & \text{if }z\in\bar\Sigma.\end {cases}
\]
If $p\in\bd\Sigma$ and local coordinates $(y_1,\dots,y_4)$ satisfy $\frac{\p r(p)}{\p y_j}=0$ for $j=1,\dots,3$, then
\[
\frac{\partial^2\tilde\delta(p)}{\partial y_j\partial y_k}=|\nabla r|^{-1}\frac{\partial^2 r(p)}{\partial y_j\partial y_k},
\]
(see for example \cite[(2.9)]{HaRa13}).  Since $T^{1,0}(b\Omega)$ is spanned by $L=\frac{\partial}{\partial z}+2i\bar z\frac{\partial}{\partial w}$, the normalized Levi-form can be computed by
\[
  \mathcal{L}(-i L\wedge\bar{L})=-|\nabla r|^{-1}=-(1+4|z|^2)^{-1/2}.
\]
Since $|L|^2=\frac{1}{2}+2|z|^2$, we have
\[
  |L|^{-2}\mathcal{L}(-i L\wedge\bar{L})=-2(1+4|z|^2)^{-3/2},
\]
and hence the only eigenvalue of the Levi-form is $-2(1+4|z|^2)^{-3/2}$.  Since this approaches $0$ as $|z|\rightarrow\infty$ for $z\in\bd\Sigma$, it follows that the Siegel upper half space is $Z(0)$ but not uniformly $Z(0)$. Thus,
(\ref{eqn:bdy dominated by dbar,dbars}) does not hold on $\Sigma$ by Theorem \ref{thm:surface int equiv to Z(q)}. This completes the proof of Theorem \ref{thm:Siegel upper half space}.\qed

% subsection: Examples
\subsection{A Positive Example}\label{subsec:examples}
%\begin{enumerate}[1.]
%\item
Let
\[
\Omega = \Big\{ (z',z_n)\subset \C^{n-1}\times\C : |z'|^2 + 2\Rre(z_n)^2 < \frac 12 \Big\}.
\]
In this case,
\[
\rho(z',z_n) = |z'|^2 + 2\Rre(z_n)^2 - \frac 12.
\]
Since $|d\rho|^2=4|z'|^2+16\Rre(z_n)^2$, we have $|d\rho|^2=8\Rre(z_n)^2+2\geq 2$ and $|d\rho|^2=-4|z'|^2+4\leq 4$ on $\bd\Omega$.  Since $i\p\dbar\rho = I$, the identity matrix, the Levi form is also an identity matrix with eigenvalues bounded between $\frac{1}{2}$ and $\frac{1}{\sqrt{2}}$.  Thus $\Omega$ satisfies $Z(q)$ uniformly for all $1 \leq q \leq n-1$.  Since $\Omega$ is also pseudoconvex, we can use the Morrey-Kohn-H\"ormander formula to prove \eqref{eqn:bdy dominated by dbar,dbars} for $t=0$.

%\item This is a variation of Example 1.
%Fixing $1 \leq q \leq n-1$ and changing $\rho$ to
%\[
%\rho(z',z_n) = \sum_{j=1}^{n-q}|z_j|^2 + \sum_{k=n-q+1}^{n-1} |z_k|^2 + 2\Big(\frac{z_n+\z_n}2\Big)^2 - \frac 12
%\]
%provides an example of domain $\Omega$ where the Levi form has $n-q$ eigenvalues of 1 and therefore satisfies $Z(q')$ uniformly for $q \leq q' \leq n-1$.

%\item
%A similar computation to Example 1 shows that an example in \cite{HaRa15u} where \textcolor{red}{Fix the next example!}
%\[
%\Omega = \Big\{ (z',z_n)\subset \C^{n-1}\times\C : |\Imm z'|^2 + 2\Rre(z_n)^2 < \frac 12 \Big\}
%\]
%also satisfies $Z(q)$ uniformly for all $1 \leq q \leq n-1$, in addition to $\dbar$ having closed range in $L^2_{0,q}(\Omega)$
%\end{enumerate}

\bibliographystyle{alpha}
\bibliography{mybib12-8-15}

\begin{thebibliography}{HPR15}

\bibitem[FK72]{FoKo72}
G.~B.\ Folland and J.~J.\ Kohn.
\newblock {\em The {Neumann} problem for the {Cauchy}-{Riemann} Complex},
  volume~75 of {\em Ann.\ of Math.\ Stud.}
\newblock Princeton University Press, Princeton, New Jersey, 1972.

\bibitem[Gan]{Gan10}
K.\ Gansberger.
\newblock On the weighted $\bar\partial$-{N}eumann problem on unbounded
  domains.
\newblock arXiv:0912.0841v1.

\bibitem[HM]{HeMcN15}
A.-K. Herbig and J.~D. McNeal.
\newblock On closed range for $\bar\partial$.
\newblock {\em submitted}.
\newblock arXiv:1410.3559v1.

\bibitem[Ho91]{Ho91}
L.-H. Ho.
\newblock {$\overline\partial$}-problem on weakly {$q$}-convex domains.
\newblock {\em Math.~Ann.}, 290(1):3--18, 1991.

\bibitem[H{\"o}r65]{Hor65}
L.~H{\"o}rmander.
\newblock ${L}^{2}$ estimates and existence theorems for the $\bar \partial $
  operator.
\newblock {\em Acta Math.}, 113:89--152, 1965.

\bibitem[HPR15]{HaPeRa15}
P.S. Harrington, M.M. Peloso, and A.S. Raich.
\newblock Regularity equivalence of the {S}zeg{\"o} projection and the complex
  {G}reen operator.
\newblock {\em Proc.~Amer.~Math.~Soc.}, 143(1):353--367, 2015.
\newblock arXiv:1305.0188.

\bibitem[HR]{HaRa15u}
P.~Harrington and A.~Raich.
\newblock Closed range for $\bar\partial$ on unbounded domains in {$\mathbb
  C^n$}.
\newblock {\em submitted}.
\newblock arXiv:1507.06211.

\bibitem[HR11]{HaRa11}
P.~Harrington and A.~Raich.
\newblock Regularity results for $\bar\partial_b$ on {CR}-manifolds of
  hypersurface type.
\newblock {\em Comm.\ Partial Differential Equations}, 36(1):134--161, 2011.

\bibitem[HR13]{HaRa13}
P.~Harrington and A.~Raich.
\newblock Defining functions for unbounded {$C^m$} domains.
\newblock {\em Rev.~Mat.~Iberoam.}, 29(4):1405--1420, 2013.

\bibitem[HR14]{HaRa14}
P.~Harrington and A.~Raich.
\newblock Sobolev spaces and elliptic theory on unbounded domains in {$\R^n$}.
\newblock {\em Adv. Differential Equations}, 19(7/8):635--692, 2014.

\bibitem[HR15]{HaRa15}
P.~Harrington and A.~Raich.
\newblock Closed range for $\bar\partial$ and $\bar\partial_b$ on bounded
  hypersurfaces in {S}tein manifolds.
\newblock {\em Ann. Inst. Fourier (Grenoble)}, 65(4):1711--1754, 2015.

\bibitem[Str10]{Str10}
E.\ Straube.
\newblock {\em Lectures on the ${\mathcal{L}}^2$-Sobolev Theory of the
  $\bar\partial$-Neumann Problem}.
\newblock ESI Lectures in Mathematics and Physics. European Mathematical
  Society (EMS), Z{\"u}rich, 2010.

\end{thebibliography}

\end{document}